\documentclass[12pt]{article}
\usepackage{amsmath, amsfonts, amssymb, amsthm, graphicx, color} 
\usepackage[utf8]{inputenc}
\usepackage[english]{babel}
\usepackage{epsf}
\usepackage{verbatim}
\newtheorem{theorem}{Theorem}[section]
\newtheorem{proposition}[theorem]{Proposition}

\newtheorem{remark}[theorem]{Remark}

\newtheorem{corollary}[theorem]{Corollary}

\renewcommand\baselinestretch{0.94}
\def\R{\mathbb R} \def\Z{\mathbb Z}
 \def\H{\mathbb H}
\def\diag{{\rm diag}\,}
\def\N{\mathbb N}
\def\C{{\mathbb C}} \def\d{{\rm d}}
\def\Q{\mathbb Q}
\def\P{\mathbb P}
\def\diam{{\rm{diam}}\,}
\def\tr{{\sf t}}
\def\norm{|\!|}
\def\supp{{{\rm supp}\,}}

\def\wt{\widetilde}
\def\Ad{\hbox{\rm Ad}}
\def\ov{\overline}
\def\Aut{\rm Aut}
\def\Aff{\rm Aff}
\def\wt{\widetilde}
\def\g{\hbox{\sf  g}}
\def\w{\hbox{\sf  w}}
\def\norm{|\! | }

\def\t{^t\!}
\def\sgn{{\rm sgn}}

\def\N{{\bf N}}

 \def\P{{\cal {P}}} 
\def\x{{\bf x}}

\def\norm{{\vert \! \vert}}

\date{}

\def\disp{\displaystyle}

\def\build#1_#2^#3{\mathrel{\mathop{\kern 0pt#1}\limits_{#2}^{#3}}}

\def\smallsquare{\vbox{\hrule\hbox{\vrule height 1 ex\kern 1
ex\vrule}\hrule}}

\def\N{{\mathbb N}}

\def\Im{\mbox{\rm Im}\,}

\renewcommand\baselinestretch{1.1}

\begin{document}

\title{On values of binary quadratic forms at integer points}
\author{ Manoj Choudhuri and S. G. Dani}
\maketitle
\begin{abstract}
We obtain estimates for the number of integral solutions in large
balls, of inequalities of the form $|Q(x,y)|<\epsilon$, where $Q$ is
an indefinite binary quadratic form, in terms of the Hurwitz continued
fraction expansions of the slopes of the lines on which $Q$ vanishes.
The method is based on a coding of  geodesics on the modular surface
via Hurwitz expansions of the endpoints of their lifts in the
Poincar\'e half-plane. 
\end{abstract}

\section{Introduction}
Consider a binary quadratic form $Q(x,y)=(ax+by)(cx+d y)$, where 
$a,b,c,d\in \R$ and $ad-bc\neq 0$. 
In this paper we exhibit a
close relationship between the growth of the number of solutions $(x,y)$, 
with $x,y\in \Z$, with $\gcd (|x|,|y|)=1$, of 

\begin{equation}
 |Q(x,y)|<\epsilon \ \ \mbox{ \rm and } \ \norm (x,y)\norm 
 \leq \rho,
\end{equation}
asymptotically as $\rho \to \infty$ and the continued fraction expansions 
of $a/b$ and $c/d$ with respect to the  Hurwitz algorithm (see below for 
details), when at least one of the ratios are irrational.  Here and in 
the sequel $\norm \cdot \norm $  stands for the Euclidean norm on~$\R^2$. 
For any set $E$ we shall denote by $\# E$ the cardinality of $E$. 

We recall (see \cite{KU} for instance) that any irrational number $\xi$ 
can be expressed as 
\begin{align*}
 \xi=a_0-\frac{\disp 1}{\disp{a_1-\frac{\disp 1}{\disp{a_2-\frac{\disp1}{\disp{\ddots}}}}}}
\end{align*}
with $a_j\in \Z$ and (i) $|a_j|\geq 2$ for all $j$ and if
(ii)  $|a_j|=2$ for some $j$ then $a_ja_{j+1}<0$ (i.e. if $|a_j|=2$ then $a_j$ 
 and $a_{j+1}$ have the opposite sign);
here the right hand side stands, as usual, for the limit (with 
existence assured) of the 
sequence of rational  
numbers represented by terms of corresponding  
truncated expressions. We draw the reader's attention to that in writing 
the expression as above the
successive terms are substracted, rather than added; this is often done in 
literature (see 
\cite{KU} for instance) on account of its being in consonance with the 
action of the modular group; an expansion as above is called 
the {\it minus continued fraction expansion} with respect to the nearest 
integer algorithm (or also as the negative
Hurwitz continued fraction expansion). We shall, as usual, denote the expression on the
right hand side by $[a_0,a_1, \dots ]$. We recall also that, 
conversely, given a sequence $\{a_j\}$ of 
integers satisfying conditions (i) and (ii) as above  there is a unique 
irrational number $\xi$
such that $\xi = [a_0,a_1, \dots ]$. 

We recall that an integral 
pair $(x,y)$ is said to be {\it primitive} if it is nonzero and $\frac 1k(x,y)$
is not an integral pair for any natural number $k$; in other words, if 
$(x,y)\neq (0,0)$ and  $\gcd (|x|,|y|)=1$. In the sequel we denote by 
$\mathfrak p$ the set of
all primitive integral pairs.  

We shall be interested in primitive solutions of inequalities as in~(1), 
with\\ $Q(x,y)=(ax+by)(cx+dy)$, where we assume the ratio $a/b$ to be irrational.
  It will be convenient to consider primitive solutions in the region 
defined by
\begin{equation}
\{ (x,y):0<|Q(x,y)|<\delta, \ cx+dy>\kappa \ 
\mbox{ \rm and } \norm (x,y) \norm \leq \rho\}, 
\end{equation}
\noindent with $\kappa >0$, as indicated  in Figure~1; estimates 
for the number of primitive integral 
solutions to inequalities as in (1) can then be obtained by putting together 
the solutions in various 
subregions (from various tentacles as in Figure~1); see Remark~1.3 for details.
 Theorem~\ref{thm:main} below addresses 
this in terms of the nearest integer continued fraction expansion of $a/b$, giving lower and upper
estimates for the number of solutions when $\rho$ is sufficiently large. 
We mention here that in the course of the proof a more concrete relationship 
is proved for a special class of quadratic 
forms called $H$-reduced forms (see \S 4 for definition of $H$-reduced, and 
Theorem~\ref{thm:reduced} for the result involved). It may also be noted that
a somewhat sharper, but more technical, version of the following theorem may
be contained in Corollary~\ref{asymptotic}; the various constants involved are 
also given in sharper form in the latter.

\begin{theorem}\label{thm:main} 
Let  $Q(x,y)=(ax+by)(cx+dy)$
be a quadratic form, where $a,b,c,d\in \R$, $ad-bc=1$, $b\neq 0$ and 
 $\displaystyle \frac ab$ 
is irrational. Let $[a_0,a_1, \dots ]$ be the continued fraction 
expansion of ${\displaystyle \frac ab}$ with respect to the nearest 
integer algorithm.  Let 
$$\alpha^-=\liminf \frac 1n \sum_{j=0}^{n-1} \log |a_j| \mbox{ \rm and } 
\alpha^+=\limsup \frac 1n \sum_{j=0}^{n-1} \log |a_j|.$$
Let $0<\delta <\frac 1\pi$ and  let $e(\delta)$ 
and $f(\delta)$ respectively denote the (asymptotic) lower density of 
$\{j\geq 0\mid |a_j|\geq \delta^{-1} +1\}$ and the upper density of 
$\{j\geq 0\mid |a_j|\geq \delta^{-1} -\frac 32\}$. Let  $\kappa >0$ be 
fixed and let 
$$D=\{(x,y)\in \R^2\mid 0<|Q(x,y)|<\delta, cx+dy> \kappa\}. $$
For any $\rho >0$ let 
$$G(\rho)=\{(x,y)\in  \mathfrak p \cap D\mid  \norm (x,y) \norm \leq \rho\}.$$
Then we have the following : 

i) if $\alpha^+<\infty$ then there exists $\rho_0$ such that 
for all $\rho \geq \rho_0$ we have  
$$\# G(\rho) \geq \frac {e(\delta)}
{(\alpha^++3)}\log \rho;$$  

ii) let $M_0= \max \{\frac 14\log \frac 95, \frac 18 \alpha^-\}$;  
then for all $M\in (0,M_0)$ and $m>f(\delta)$,
there exists $\rho_0$ such that for all $\rho \geq \rho_0$ we have
$$\# G(\rho)\leq \frac mM \log \rho.$$ 
(if $\alpha^- $ is infinite, the assertion holds for all positive $M$).

\end{theorem}

\begin{figure}
 \centering
 \includegraphics[scale=0.8]{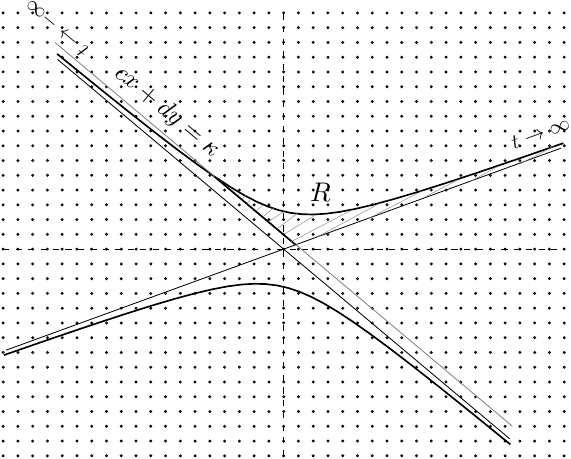}
 \caption{The region $R=D\cap \{v\in\R^2\mid Q(v)>0\}$}
 \end{figure}

\begin{remark}
{\rm The assertion in (ii) shows in particular that  
 for any quadratic form $Q$ as in Theorem~\ref{thm:main} 
and $0<\delta <\frac 1\pi$, such that for all sufficiently large $\rho$ we have 
$\# G(\rho)\leq (\frac 14\log \frac 95) \log \rho$. The area of the of the 
region $D\cap B(\rho)$, where $B(\rho) =\{(x,y)\mid \norm (x,y)\norm \leq \rho\}$,
$D$ is as in the theorem and for some choice of 
$\kappa >0$, is 
asymptotic to $\log \rho$, and the preceding observation signifies that the 
the number of primitive integral pairs in $D\cap B(\rho)$ is bounded above, for all 
large $\rho$, by a fixed multiple of the area of the region. We note that this 
contrasts the situation for linear forms in place of the quadratic forms 
where the number of primitive integral points in analogously defined regions 
can be arbitrarily large in proportion to the area of the region, depending 
on the linear form, on account 
of existence of very well approximable numbers. We may also mention 
here the following, concerning quadratic forms in higher number of variables. 
  For all nondegenerate
indefinite,  irrational (not a scalar multiple of a
form with rational coefficients)  forms in $n\geq 3 $ variables 
such a ratio is bounded below and for $n\geq 5$ the ratio converges to~$1$; 
for  $n=3$ or $4$ there are counter-examples to the second part; 
see~\cite{M} for details. } 
\end{remark}

\begin{remark}
{\rm 
Let $Q$ be the quadratic form as in Theorem~\ref{thm:main} such  
that both $a/b$ and $c/d$ are irrational. Let 
$0<\delta <\displaystyle \frac 1\pi$ and for any $\rho >0$ let 
$$~~~~~H(\rho)=\{(x,y)\in \mathfrak p\mid 0<|Q(x,y)|<\delta
\mbox{ \rm and } \norm (x,y) \norm \leq \rho\}, \mbox{\rm and}$$\  
$$G'(\rho)=\{(x,y)\in \mathfrak p\mid 0<|Q(x,y)|<\delta, ax+by>\sqrt \delta
\mbox{ \rm and } \norm (x,y) \norm \leq \rho\}.$$
Then for any $\rho >0$ we have that $\# H(\rho)$ differs 
from $2(\# G(\rho)+\# G'(\rho))$ by at most~$2$, as maybe observed from 
Figure~1; note that here 
$\kappa$ is chosen to be $\sqrt \delta$. Also, 
$\# G'(\rho)$ admits analogous estimates as $\# G(\rho)$ 
with the continued fraction expansion of $c/d$ in place of $a/b$. Thus 
the theorem provides estimates for $\# H(\rho)$ under appropriate 
conditions as in the hypothesis of the theorem. 

}
\end{remark}

\begin{remark}\label{remark 1.4}
{\rm It may be emphasized that Theorem~\ref{thm:main} applies to {\it all} 
irrational numbers $\alpha$; thus, given any sequence $\{a_j\}$ with 
$a_j\in \Z$ such that 
$|a_j|\geq 2$ for all $j$ and $a_ja_{j+1}<0$ if $|a_j|=2$ we have an 
 irrational number $\alpha=[a_0,a_1, \dots ]$ and a corresponding 
result for any $c,d \in \R$ such that $ad-bc=1$. One may ask what  the
typical, or generic, values of the constants involved are.
Specifically the question may be formulated as follows. 
As in the case of the usual continued fractions there is a corresponding 
Gauss map associated to the negative Hurwitz continued fraction, defined on $I=[-\frac 12, 
\frac 12]\backslash \Q$, by $T(x)=-\frac 1x +\nu (x)$, where $\nu (x)$ denotes 
the integer nearest to $x$. An invariant probability measure $\mu$ for the map, equivalent
to the Lebesgue measure is described in \cite{NIT}. It is however not in a closed form.
The measure $\mu$ can also
be obtained using the representation of the geodesic flow as a flow built under a 
function, with respect to a cross-section as in \cite{KU2}.
It follows from the ergodicity of the geodesic flow associated to the modular 
surface that the measure  $\mu$ is ergodic with respect to $T$.
The ergodicity implies that for almost all $\alpha =[a_0,a_1, \dots ]$, 
$\frac 1n \sum_{j=0}^{n-1} \log |a_j|$ converges to $\int f d\mu$, where 
$f:I\to (0,\infty)$ is the function defined by 
$f(x)=\log  |a|$ if  
$x\in \left (\frac 1{a+\frac12}, \frac 1{a-\frac12}\right )$.
(see the proof of (3.26) in \cite{EW}); 
we note that $f$ is integrable and $\int f d\mu$ is a positive constant; 
this constant then is the generic value of 
$\alpha^+$ and $\alpha^-$ as in the statement of the theorem. Similarly, for 
a $\delta >0$ the generic values of $e(\delta)$ and $f(\delta)$ are seen to 
be $\mu([-\frac 1k,\frac 1k])$ and $\mu([-\frac 1l,\frac 1l])$ respectively, 
where $k=[\delta^{-1}+1]$ and $l=[\delta^{-1}-\frac 32]$.  

}
\end{remark}
\bigskip
\begin{center}
 {\bf Summary of revisions over the previous version}
\end{center}

The following changes have been made over the previous version :
\begin{itemize}
 \item [$\langle 1 \rangle$] In the definition of $e(\delta)$ and $f(\delta)$, $2\delta^{-1}$
 is replaced by $\delta^{-1}$, the same changes are involved in giving generic 
 values of $e(\delta)$ and $f(\delta)$ in Remark \ref{remark 1.4}. 

 \item[$\langle2\rangle$] Statement $(ii)$ of Theorem \ref{thm:main} has been modified to  
 achieve greater clarity. We note that 
the statement in the review of the present paper in Math. Reviews, by David Simmons, that
``one of the inequalities in (ii) should be flipped:
$M\geq \frac 14\log\frac 95$ should have the opposite inequality." is not correct; the original assertion is justified,
by Corollary~5.1, and is reinterpreted in the statement of the theorem above. By the same Corollary 
$M=\frac 18 \alpha^-$ is also allowable when $\alpha^-<\infty$ (this is not included in the main statement since 
$\alpha^-$ can be infinite). When $\frac{1}{4}\log\frac{9}{5}
<\frac{1}{8}\alpha^-$, taking $M=\frac{1}{8}\alpha^-$ or any value between $\frac{1}{4}\log\frac{9}{5}$ and 
$\frac{1}{8}\alpha^-$ yields a better estimate than with $M\leq\frac{1}{4}
\log\frac{9}{5}$; thus if one restricts to $M\leq\frac{1}{4}\log\frac{9}{5}$ as suggested in the Review, the statement 
would not cover the stated result.
\item[$\langle 3 \rangle$] The description of the invariant measure for the Gauss map given in Remark \ref{remark 1.4}
of the previous version turns out to be incorrect. The remark has now been modified and a reference is given.
\end{itemize}

\bigskip
\noindent{\it Acknowledgement}: The authors are thankful to the referee 
for suggestions leading to improvement of the presentation of the paper.
The authors are also thankful to  L.~Singhal for his help in drawing the 
figures included in this paper.

The authors are thankful to David simmons for pointing out the errors 
noted in items $\langle 1 \rangle$ and $\langle 3 \rangle$ in the summary noted above.
\section{A correspondence}
Let $G=SL(2,\R)$. 
We denote by $\{e_1,e_2\}$ the standard basis of $\R^2$. 
Let $Q_0$ be the quadratic form on $\R^2$ defined by $Q_0(xe_1+ye_2)= xy$ for all
$x,y\in \R$. For $g\in G$ we denote by $Q_g$ the quadratic form defined by 
$Q_g(v)=Q_0(g^{-1}v)$ for all $v\in \R^2$.

For $t\in \R$ we denote by $a_t$ the matrix $\diag (e^{t/2}, e^{-t/2})
\in G$.  For $g\in G$ and $t\in \R$, if 
$g_t= ga_tg^{-1}$, then for any $v\in
\R^2$ we have $Q_g(g_t v)=Q_0(g^{-1}g_tv)=Q_0(a_tg^{-1}v)=Q_0(g^{-1}v)
=Q_g(v)$; thus $\{g_t\}_{t\in \R}$ is contained in $SO(Q)$, and in fact 
coincides, by dimension considerations, with $SO(Q)^+$, the connected 
component of the identity in $SO(Q)$. 
In particular $\{v\in \R^2\mid Q_g(v)\neq 0\}$ is $\{g_t\}$-invariant.
We note that 
 in each connected component of $\{v\in \R^2\mid Q_g(v)\neq 0\}$ the orbits 
of $\{g_t\}_{t\in \R}$ are the level curves of $Q_g$ and they are asymptotic 
to the pair of lines defined by $Q_g(v)=0$. We denote by  
$L_g^+$ and  $L_g^-$ the linear forms on $\R^2$ such that 
$Q_g=L^+_gL^-_g$ and the level curves of $Q_g$, viewed as orbits of $\{g_t\}$, 
are asymptotic to $L_g^+(v)=0$ as $t\to \infty$ and  $L_g^-(v)=0$ as 
as $t\to -\infty$; see Figure~1.

Through the rest of the section we fix a $g\in G$ and let $Q=Q_g$. It 
may be mentioned that the results will have nontrivial content only when 
at least one of $L_g^+$ and  $L_g^-$ is not a rational form, but we make
no specific assumption in this respect.  Let $g_t=ga_tg^{-1}$ for all 
$t\in \R$. For any $v\in \R^2$ and any subset $C$ of 
$\R^2$ we define $$ I_v(C)= \{t\geq 0 \mid v\in g_tC\}\mbox{ \rm and }
R(C)=\bigcup_{p\in \Z^2\backslash \{0\}}I_p(C).$$ 
\begin{proposition}\label{prop0}
Let $C$ be a convex subset of $\R^2$ containing $0$ and with area less than 
$\frac 12$. Then the following holds: 

i) for any $v\in \R^2$, $I_v(C)$ is an interval in $\R$;

ii) for $p,p'\in \Z^2$ ,  $I_p(C)$ and $I_{p'}(C)$ are 
contained in disjoint connected components of $R(C)$ if and only if 
$p$ and $p'$ are linearly independent; 

iii) there exists $\kappa >0$ such that the distance 
between any two (successive) connected components of $R(C)$ is at 
least $\kappa$. 

\end{proposition}

\proof i) If $g_tv, g_{t'}v\in C$ for $t< t'\in \R$, then 
$\{g_sv\mid t\leq s\leq t'\}$ is a segment of a hyperbola, and when $C$ is
a convex set containing $0$ the segment is contained in $C$. This proves (i).

\noindent ii) Now let $p,p'\in \Z^2$ and 
suppose $I_p(C)$ and $I_{p'}(C)$ are contained in the same connected component
of $R(C)$.  Let 
$t\in I_p(C)$ and $t'\in I_{p'}(C)$ be given, with say $t<t'$; then $t$ and 
$t'$ are contained in the same connected component of $R(C)$. Since $R(C)$ 
is the union of intervals of the form $I_q(C)$, $q\in \Z^2\backslash \{0\}$, 
$t$ and $t'$ being contained in the same connected component implies that 
there exist  $p=p_0, p_1, \dots ,
p_{k-1},p_k=p' \in \Z^2$ and $t=t_0 <t_1<t_2 \dots <t_k=t'$ such that $p_j\in 
g_{s}C$ for all $s\in [t_{j}, t_{j+1}]$ and  $j=0, \dots, k-1$. 
Consider any $1\leq i\leq k$. 
We have $g_{t_{j}}^{-1}p_{j-1}, g_{t_{j}}^{-1}p_{j} \in C$, and
hence by the condition in the hypothesis the triangle with 
vertices at $g_{t_{j}}^{-1}p_{j-1}, g_{t_{j}}^{-1}p_{j}$ and $0$ 
has  area less than  
$\frac 12$. It follows that the 
triangle with vertices as $p_{j-1},p_j$ and $0$ has area less that $\frac 12$. 
Since these vertices are integral points the conclusion implies that 
$p_j$ and $p_{j-1} $ are linearly dependent, for all $j=1,\dots, k$. 
In particular we get that $p$ and $p'$ are linearly dependent. 

Now suppose $p,p'$ are linearly dependent,
say $p=kp'$ with $k>1$. Let $t\in I_p(C)$. Then there exists 
$v\in C$ such that $p=g_t(v)$. Then 
$p'=k^{-1}p=k^{-1}g_t(v)\in g_tC$, since $k^{-1}<1$. 
This shows that $I_p$ is contained in $I_{p'}$. In particular 
$I_p(C)$ and $I_{p'}(C)$ are contained in the same connected component of 
$R(C)$.  This proves~(ii). 

\noindent iii) We can find a convex neighbourhood $C'$ of $C$ with area less 
than $\frac 12$. Then there exists a $\kappa >0$ such that $g_tv\in C'$ for 
all $v\in C$ and $-\kappa \leq t \leq \kappa$. We see that each connected 
component of $R(C)$ are contained in a unique connected component of $R(C')$. 
Therefore the successive connected components of $R(C)$ are at a distance
bounded below by $\kappa$. This proves (iii) \qed

\begin{proposition}\label{prop1}
Let $Q=Q_g$ and $\theta> 0$ and $0\leq r <r',$ be given. Let 
$$ \Omega=\{ v\in
\R^2\mid 0 < L_g^+ (v) < \theta L_g^-(v) \mbox{ \rm and }r<Q(v) <r'\}
\mbox{ \rm and } S=\Omega \cap \Z^2.$$ 
Let $\sigma \geq \sqrt {r'/\theta}$, and let  
$$T= \{v\in \R^2 \mid  0 < L_g^+ (v) < \theta L_g^-(v) 
\mbox{ \rm and }  0\leq  L_g^-(v)\leq \sigma\}.$$
Suppose that $T$ has  area less than  $\frac 12$. 
 For $p\in S$  let $I_p=\{t\geq 0\mid p\in g_tT\}$.  Then  
 each $I_p(T)$, $p\in S$, is  nonempty, and  for $\tau > \sigma $ 
any maximal set of pairwise linearly 
independent vectors contained in 
$\{p\in S\mid L_g^-(p)\leq \tau\}$ has cardinality equal to the 
number of connected components of  $[0,2\log (\tau/\sigma)]\cap 
R(T\cap \Omega)$. 
\end{proposition}

\proof The assumption $\sigma \geq \sqrt {r'/\theta}$ ensures that
the set of values of $Q$ on the triangle $T$ 
contain the interval $(0, r')$. This implies that for $p\in S$ 
there exists $v\in T$  and $t\geq 0$ such that $p=g_tv$, 
showing that  $I_p(T)$ is nonempty. 

Now, for  $\tau>0$  let  $S(\tau) =\{ p\in S\mid L_g^-(p)
\leq \tau\}$. 
Now consider $\tau >\sigma$ and $p\in S(\tau)$. Then there exists 
$v\in T\cap \Omega$ with $L_g^-(v)=\sigma$ and $t>0$ such that $p=g_t(v)$.
Thus  $L_g^-(p)=e^{t/2}L_g^-(v)=e^{t/2}\sigma$ and hence $t\leq 
2\log L_g^-(p)/\sigma\leq 2\log (\tau/\sigma) $. Therefore for each 
$p\in S(\tau)$ there exists a $t\in [0,2\log (\tau/\sigma )]$ such that 
$p\in g_t(T\cap \Omega)$; let $J_p$ denote the connected component of 
$ [0, 2 \log (\tau /\sigma)]\cap R(T\cap \Omega)$ containing $t$. 
By Proposition~\ref{prop0} 
if $p,p'\in S(\tau)$ are linearly independent then $I_p(T)$ and $I_{p'}(T)$ 
belong to disjoint connected components of $R(T)$ and hence $J_p$ 
and $J_{p'}$ are disjoint. This 
shows that the number of elements in any set of linearly independent 
vectors in $S(\tau)$ is bounded by the number of connected components of
$ [0, 2  \log (\tau /\sigma)]\cap R(T\cap \Omega)$. 

Now let $J$ be any connected component of $ [0, 2 
\log (\tau /\sigma)]\cap R(T\cap \Omega) $, and  
 let $t\in J$. Then there exists $p\in \Z^2$ and $v\in T\cap \Omega$ such that 
$p=g_t(v)$. Clearly there exists $t'\in J$ and $v'\in T\cap \Omega$ such that 
$L_g^-(v')=\sigma $ and $g_t(v)=g_{t'}(v')$, and hence by modifying 
notation we may assume that $L_g^-(v)=\sigma $. 
Hence  $L_g^-(p)=e^{t/2}L_g^-(v)= e^{t/2}\sigma$. 
As $t\leq 2\log (\tau /\sigma) $, we get that 
$L_g^-(p)\leq \tau $ and hence  $p\in S(\tau)$. Now suppose $t$ 
and $t'$ belong to different connected components, and let $p,p'
\in S(\tau) $ be the elements obtained as above corresponding to
$t$ and $t'$ respectively. Then $t\in I_p$ and $t'\in I_{p'}$. 
Thus $I_p(T)$ and $I_{p'}(T)$ are intervals containing respectively 
$t$ and $t'$ belonging to distinct connected components of 
$R(T\cap \Omega)$ and hence by Proposition~\ref{prop0} 
$p$ and $p'$ are linearly independent. This shows that  
there are at least as many  mutually linearly independent 
vectors in $S(\tau)$ as the number of connected components
of $[0, 2 \log (\tau /\sigma)]\cap R(T\cap \Omega)$, which proves 
the proposition. \qed 

\medskip
If $p$ is a primitive integral vector, then $p=xe_1+ye_2\in \Z^2$ with  gcd $(|x|,|y|)=1$. 
For any $\sigma >0$ let $$W(\sigma)=\{v=xe_1+ye_2 \in \R^2\mid 0< y
< x\leq \sigma\}.$$
For $\tau >\sigma >0$ we denote by $n(\tau, \sigma)$ the the number of 
connected 
components of $\{t\in  [0,2\log (\tau/\sigma) ]\mid ga_t W(\sigma)\cap 
\Z^2\neq \emptyset \}$ or equivalently of the set of $t\in [0,2\log (\tau/\sigma)]$ 
such that $ a_{-t}\lambda \in W(\sigma)$, for some $\lambda \in \Lambda$, where 
$\Lambda :=g^{-1}\Z^2$. 

\begin{figure}[h]
\centering

\includegraphics[scale=0.8]{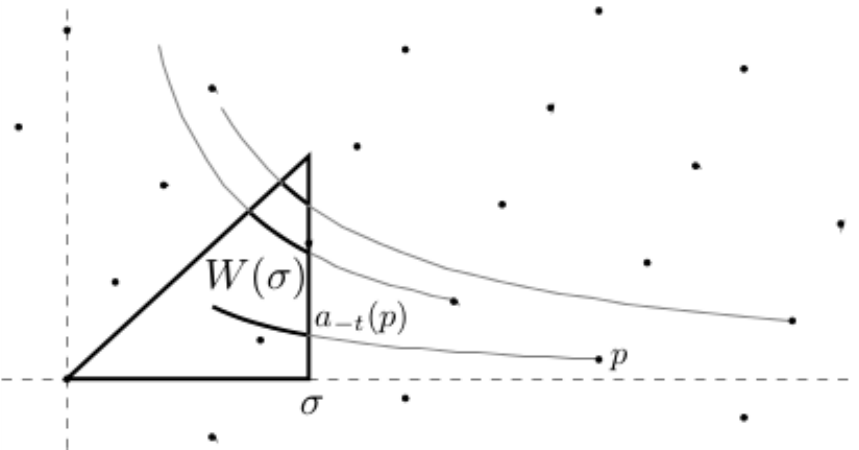}

\caption{Trajectories passing through $W(\sigma)$}

\end{figure}

\begin{corollary}\label{cor1}
Let  $Q=Q_g$.  Let  $0<\epsilon <1$ be given. For $\rho >0$ let 
$$F(\rho) =\{p\in \mathfrak p\mid  0<L_g^+(p) < L_g^-(p), \ 0<Q(p)<\epsilon, \ 
\norm p \norm \leq \rho\}.$$ Then for all sufficiently large $\rho $ 
the cardinality of $F(\rho)$ differs from  $n(\norm ge_1\norm ^{-1} \rho, \sqrt 
\epsilon)$ by at most~$1$. 

\end{corollary}

\proof In Proposition~\ref{prop1} we choose $\theta =1$, $r=0, 
r'= \epsilon$, and let $\Omega$ and $T$ be as in that Proposition with $\sigma=\sqrt{r'}$.  
Then for the choices as above  $\Omega$ contains $T$. Also,  since
 $\epsilon <1$  the area  of $T$ is less than $\frac 12$. We note also
that in the above notation $T$ coincides with $gW(\sqrt \epsilon)$, and hence 
$g_tT=ga_tW(\sqrt \epsilon)$. Then  by
Proposition~\ref{prop1}, for any $\tau >\sqrt \epsilon$ the number of 
$p$ in $\mathfrak p$ 
such that $ 0<L_g^+(p) < L_g^-(p) \leq \tau$ and $0<Q(p)<\epsilon$, is
$n(\tau, \sqrt \epsilon)$. Let $\delta >0$ be arbitrary. We note that 
for $v\in T$ and $t$ sufficiently large 
$$|\,\norm g(e_1)\norm ^{-1}\norm g_t(v)\norm - L_g^-(g_tv)\,|<\delta. $$
Hence we get that for all large $\rho$ the cardinality of  $F(\rho)$ 
is bounded between \\ $n(\norm g(e_1)\norm ^{-1}\rho -\delta,  \sqrt \epsilon)$ 
and  $n(\norm g(e_1)\norm ^{-1}\rho +\delta, \sqrt \epsilon )$. By 
Proposition~\ref{prop0} the successive connected components are at 
distance bounded below by
a positive constant, and hence   when $\delta $ is sufficiently
small, for all large $\rho$,  $n(\norm g(e_1)\norm ^{-1}\rho +\delta, 
\sqrt \epsilon )$ is at most one more than $n(\norm g(e_1)\norm ^{-1}
\rho -\delta, \sqrt \epsilon )$.  Hence  the cardinality of 
$F(\rho)$ differs from  $n(\norm g(e_1)\norm ^{-1}\rho,\sqrt 
\epsilon )$ by at most $1$. This proves the corollary. \qed

For $\delta >0$ let
 $B_\delta$ be the open ball of radius $\delta$ in $\R^2$, centered at $0$, 
(in the usual metric) and let $B_\delta'=\{(x,y)\in B_\delta \mid 0 <y<x\}$. 

\begin{corollary}\label{cor2}
Let  $Q=Q_g$  and $0<\delta <\frac 2{\sqrt \pi} $ 
be given. Then for all sufficiently large  $\rho >0$ the cardinality of 
$$
\{p\in \mathfrak p\mid  0<L_g^+(p) \leq L_g^-(p), \ 0<Q(p)<\frac 12 \delta^2, \ 
\norm p \norm \leq \rho\}$$  differs from  
the  number of connected components of $\{t\in [0, 2\log (\sqrt 2 \delta^{-1}
\norm ge_1 \norm ^{-1} \rho) ] \mid ga_t B_\delta'\cap \Z^2\neq \emptyset\}$ by at most $2$. 
\end{corollary}

\proof 
In view of  Corollary~\ref{cor1}, for $\delta <1$  the cardinality of the set
under consideration differs 
by at most one from the  number of connected components of 
$\{t\in [0, 2\log  (\sqrt 2\delta^{-1}\norm ge_1\norm ^{-1} \rho  ]\mid ga_t 
W(\delta/\sqrt 2)\cap \Z^2\}$. 
Now, each of these connected components is contained in a connected 
component of $\{t\in [0, 2\log (\sqrt 2 \delta^{-1}\norm ge_1\norm ^{-1} \rho) ]
\mid ga_t B_{\delta}' \cap \Z^2\neq \emptyset \}$. Moreover when $\delta <
\frac 2{\sqrt \pi}$ the area of $B_{\delta}'$  is less than $\frac 12$, 
and hence by Proposition~\ref{prop0} distinct connected 
components of the former are contained in distinct connected 
components of the latter. Also, at most one connected component
of  $\{t\in [0, 2\log (\sqrt 2\delta^{-1}\norm ge_1\norm ^{-1} \rho)  ]
\mid ga_t B_{\delta}' \cap \Z^2\neq \emptyset\}$ can fail to intersect 
$\{t\in [0, 2\log (\sqrt 2 \delta^{-1}\norm ge_1\norm ^{-1} \rho ) ]
\mid ga_t W(\delta/\sqrt 2)\cap \Z^2\neq \emptyset\}$.
The assertion as in the corollary is now immediate from these observations. \qed

\medskip
We shall say that two function $f$ and $f'$ on $(0,\infty)$
are {\it comparable} if the function $|f-f'|$ is bounded over $(0,\infty)$. 

\begin{corollary}\label{cor3}
 Let  $Q=Q_g$  and $0<\delta <\sqrt \frac 2\pi $ and $\kappa>0$
be given. For $\tau>0$ let $c^+_g(\tau)$ and $c_g(\tau)$ denote, respectively, 
the number of connected components of $\{t\in [0, \tau ]\mid ga_t 
B_\delta\cap \Z^2\neq \emptyset\}$ and  $\{t\in [-\tau , \tau ]\mid ga_t 
B_\delta\cap \Z^2\neq \emptyset\}$. 

i) for any 
$\rho >0$  the cardinality of $$G^+(\rho)= \{p\in \mathfrak p\mid 
 \ 0<|Q(p)|<\frac 12 \delta^2,  L_g^-(p)>\kappa, \ \norm p \norm \leq \rho\}$$  
is comparable to  $c^+_g (2\log \rho )$.  

ii) for any $\rho>0$
the cardinality of $$G(\rho)= \{p\in \mathfrak p\mid 
 \ 0<|Q(p)|<\frac 12\delta^2, \ \norm p \norm \leq \rho\}$$ is comparable to 
$2c_g (2\log \rho )$. 

\end{corollary}

\proof i) 
The intersection of $G^+(\rho)$ with $\{v\in \R^2\mid L^+_g(v)>0\}$ differs
from the set as in Corollary~\ref{cor2} by only a finite set (consisting 
of elements of $\mathfrak p$ contained in the compact subset $\{v\in \R^2\mid 
0 \leq L^+_g(v) \leq L^-_g(v)\leq \kappa \mbox{ \rm and } Q(v)\leq\frac 12 
\delta^2\}$)
and hence by Corollary~\ref{cor2} the cardinality 
of $G^+(\rho)$ is comparable to  the  number of connected components of 
$\{t\in [0, 2\log  (\sqrt 2 \delta^{-1} \norm ge_1\norm ^{-1} \rho) ]\mid ga_t 
B_\delta'\cap \Z^2\neq \emptyset\}$. Since the distance between successive 
connected components is bounded below (by Proposition~\ref{prop0}), the 
latter is comparable to   the  number of connected components of 
$\{t\in [0, 2\log \rho ]\mid ga_t B_\delta'\cap \Z^2\neq \emptyset\}$.
Similarly we obtain that the cardinality of $G^+(\rho) \cap \{v\in \R^2\mid 
L^+_g(v)<0\}$ is comparable to the number of connected components of 
$\{t\in [0, 2\log \rho ]\mid ga_t B_\delta''\cap \Z^2\neq \emptyset\}$, 
where $B_\delta''=\{(x,y)\in B_\delta \mid 0 <-y<x\}$. We note  
that since 
$\delta < \sqrt {\frac 2 \pi}$ the area of the convex closure of $B_\delta'
\cup B_\delta''$ is less than $\frac 12$ and hence the connected components
of the sets $\{t\in [0, 2\log \rho ]\mid ga_t B_\delta'\cap \Z^2\neq \emptyset\}$
and $\{t\in [0, 2\log \rho ]\mid ga_t B_\delta''\cap \Z^2\neq \emptyset\}$
are disjoint from each other. Also, their union is the set of connected 
components of $\{t\in [0, 2\log \rho ]\mid ga_t B_\delta\cap \Z^2\neq 
\emptyset\}$, except possibly for one connected component in the latter
corresponding to a $p\in \mathfrak p$ such that $L_g^-(p)=0$, in case $L_g^-$ is 
a rational linear form. Hence the cardinality of $G^+(\rho)$ is asymptotic 
to $c_g^+(2\log \rho)$. This proves assertion~(i). 

\noindent ii) Since when $p\in \mathfrak p$ 
belongs to $G(\rho)$  so does 
$-p$, from~(i) we get also that the cardinality of 
$\{p\in \mathfrak p\mid 
 \ 0<|Q(p)|<\frac 12 \delta^2,  |L_g^-(p)|>\kappa, \ \norm p \norm 
\leq \rho\}$ is $2 c_g^+(2\log \rho)$. 

Analogously, the set  $\{p\in \mathfrak p\mid 
 \ 0<|Q(p)|<\frac 12 \delta^2,  |L_g^+(p)|>\kappa, \ \norm p \norm \leq \rho\}$ 
has cardinality comparable to the number of connected components of 
$\{t\in [- 2\log \rho ,0]\mid ga_t B_\delta\cap \Z^2\neq \emptyset\}$. 
We note that $G(\rho)$ differs from $G^+(\rho)\cup G^-(\rho)$ 
by a fixed finite set. 
Therefore we get that the cardinality of $G(\rho)$ is comparable to  
$2c_g(2\log \rho)$. 
This proves the corollary. \qed

\section{The geodesic flow and continued fraction expansions}

We next relate the conclusion in Corollary~\ref{cor2} to the geodesic flow 
associated to the modular surface. The reader is referred to \cite{KU} for 
various general results on this topic used in the sequel. 

Let $\H^2=\{z\in \C\mid \Im
z>0\}$ be the Poincar\'e upper half plane. We view $\R \cup \{\infty\}$ as
the boundary $\partial \H^2$ of $\H^2$. We consider 
$\H\cup \partial \H^2$ equipped with the usual action of $G= SL(2,\R)$. 
We recall that the geodesics $\{\varphi_t\}$ in $\H^2$ are Euclidean 
semicircles joining a pair of (distinct) points in $\R \cup  \{\infty\}$, 
the boundary of 
$\H^2$; the two end points, say $u$ and $w$, are the limits of $\{\varphi_t\}$
as $t\to -\infty$ and 
$t\to \infty$, and are respectively called the {\it repelling} and 
{\it attracting end points} of the geodesic; for simplicity we may 
identify the geodesic as the geodesic joining $u$ and $w$.

Let $K$ be the subgroup of $G$ consisting of elements fixing $i$ under 
the action on $\H^2$, namely the elements of $G$  acting as rotations on 
$\R^2$. The $G$-action on $\H^2$ is transitive, and as such $\H^2$ can 
be realised as $G/K$. The geodesics 
in $\H^2$ then correspond to $\{ga_tK\}_{t\in \R}$, $g\in G$ (see \cite{K}, 
for instance), with the latter as the geodesic joining  $g(0)$ and 
$g(\infty)$.

We denote by $N$ 
the subgroup of $G$ consisting of all upper triangular unipotent matrices. 
The orbits of the $N$-action on $\H^2$ consist of  
horizontal lines. We shall also use the notation  
$A=\{a_t\mid t\in \R\}$ and, for $\delta >0$, $A_\delta= \{a_t\mid {t}< 
\log \delta \}$. 

 In the context of Corollary~\ref{cor3} we would be interested, 
for given $g\in G$, in solutions 
of $a_t^{-1}g^{-1}p \in B_\delta$ with $p\in \mathfrak p$ and $t\in \R$. Let $\Gamma
=SL(2,\Z)$. Then $\mathfrak p=\{\gamma e_1\mid 
\gamma \in \Gamma\}$. Let $\gamma \in \Gamma$ be such that $p=\gamma e_1$. 
Then the condition as above translates to $a_t^{-1}g^{-1}\gamma e_1
\in B_\delta$ which is equivalent to  $a_t^{-1}g^{-1}\gamma \in KA_\delta N$
and in turn to $\gamma^{-1}ga_t \in NA_\delta^{-1}K$, upon taking inverses. 
For $\delta >0$ let $\H_\delta=\{x+iy\mid y>\delta^{-2}\}$. Then 
considering the $G$-action on $\H^2$ we see that the condition 
as above is equivalent to  $\gamma^{-1}ga_t (i)\in \H_\delta$. Now consider 
the quotient map, say $\eta$, of $\H^2$ onto the ``modular surface'' 
$M=\Gamma \backslash \H^2$ and for $\delta >0$ let $M_\delta =\eta(\H_\delta)$. 
Then the above condition is equivalent to $\eta (ga_t(i))\in M_\delta$. 
Thus we 
see that 
$$\{t\in \R\mid a_t^{-1}g^{-1}p\in B_\delta \mbox{ \rm for some }p\in \mathfrak p\}=
\{t\in \R\mid \eta (ga_t(i))\in M_\delta \}. \quad {(*)}$$

Now $ga_t(i)$ is a geodesic in $\H^2$ and $\eta (ga_t(i))$ is its image 
under the quotient map (it is a geodesic with respect to the induced metric
on $M$, but we will not need to go into the geometry on the quotient). The 
above observation 
 enables us, using  Corollary~\ref{cor2}
on the one hand and coding of geodesics on the other hand 
to count the number of primitive solutions of quadratic inequalities
in large balls in $\R^2$. 

Before proceeding with the main results 
we note the following: 

\begin{remark}\label{rem1}
{\rm 
If $g\in G$ and $Q=Q_g$, and if the inequality 
$|Q(p)|<\epsilon $ admits solutions $p\in \Z^2$ for all $\epsilon>0$ then 
sets on the left hand side of $(*)$ have to be nonempty for all $\delta>0$ (we  
shall not concern ourselves with the precise correspondence between the 
values here). Hence the 
observation shows in particular that if 
the image of the geodesic joining $g(0)$ and 
$g(\infty)$ under the quotient map onto the modular surface $M=\Gamma 
\backslash \H^2$ is bounded in $M$,  then there exists 
$\epsilon >0$ such that $|Q(p)|< \epsilon$ has no nonzero solution $p$ in 
$\Z^2$.}
\end{remark} 

We now begin by introducing some notation, and recalling some definitions and 
facts about geodesics in $\H^2$ and their images in $M$. 

Let $T^1(\H^2)$ be the unit tangent bundle of $\H^2$, viewed as the set of 
pairs $(z,\zeta)$ with $z\in \H^2$ and $\zeta $ a unit tangent direction 
at $z$. For $(z,\zeta) \in T^1(\H^2)$ we denote by 
 $\tilde \varphi (z,\zeta)$ 
the geodesic $\varphi= \{\varphi_t\}$ such that $\varphi_0=z$ and 
$\zeta$ is the tangent direction to $\varphi$ at $z$.

We recall, from Katok and Ugarcovici~\cite{KU}, that a geodesic $\varphi 
=\{\varphi_t\}$  in $\H^2$ 
with $u$ and $w$ the corresponding repelling and attracting endpoints in 
$\R\cup \{\infty\}$  is said to be {\it $H$-reduced} if 
$|w|>2$ and $\sgn \,(w)u\in[\lambda-1, \lambda ]$, where $\lambda 
=\frac 12(3-\sqrt{5})$; in \cite{KU} the notation is $r$ in place of the 
$\lambda$ used here. 
Now let $\mu = (23-3\sqrt 5)/22$ and let $C$ be the arc in $\H^2$ defined 
by $$C= \{z=x+iy \in \H^2\mid |z|=1, |x|<\mu\}.$$ 
We note that $\mu$ chosen here is the $x$-coordinate of the point of 
intersection 
of the geodesic joining $\lambda$ and $2$ with $\{z\in \H^2\mid |z|=1\}$. 
It is straightforward to 
verify that for every $H$-reduced geodesic $\varphi=\{\varphi_t\}$ there 
exists a (unique) $t_0\in \R$ such that $\varphi_{t_0}\in C$.  

\begin{figure}
 
\centering
 \includegraphics[scale=1]{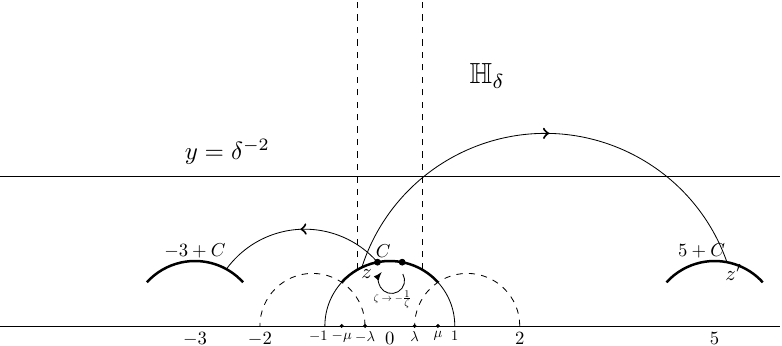}
 \caption{Segmentation of the geodesic trajectories}
 \end{figure}

We shall be interested in the set of tangent vectors defined by
$$\Phi =\{(z,\zeta)\in T^1(\H^2)\mid z\in C, \tilde \varphi (z,\zeta) 
\mbox{ \rm is } H\mbox{-\rm reduced} \}.$$
We recall that every geodesic in $\H^2$ is equivalent, under the 
$\Gamma$-action, to a $H$-reduced geodesic~\cite{KU}.
Thus every geodesic in $\H^2$ is equivalent under 
the $\Gamma$ action to (and hence has the same image in $M$ as)
a geodesic of the form $\tilde \varphi (z,\zeta)$, $(z,\zeta) \in \Phi$. 
It is known  that $\Phi$ as above corresponds to a 
cross-section for the geodesic flow  (see \cite{KU}), though 
it is not stated 
in this form. We shall first give here a direct verification of this, 
together information on the return times that we will be using
below.

Let $(z_0,\zeta_0) \in \Phi$ and $\varphi =\tilde \varphi (z_0,\zeta_0)=
\{\varphi_t\}$ be the corresponding geodesic. Let $u$ and $w$ be the
repelling and attracting endpoints of $\varphi$. We shall suppose that 
$w$ is an irrational number. We now consider the (one sided) trajectory 
$\{\varphi_t\}_{t\geq 0}$, and divide it into segments as follows. 
Let $[a_0,a_1,\dots ] $ be the 
minus continued fraction expansion of $w$ following the Hurwitz 
algorithm.

We claim that the trajectory $\{\varphi_t\}_{t\geq 0}$
intersects the arc $\{a_0+z\mid z\in C\}$, the translate of $C$ at $a_0$. 
Suppose first that $a_0>0$, so $a_0\geq 2$; see Figure~3, for reference
for the following argument, where $a_0$ is chosen to be $5$. Since $w\in (a_0-
\frac 12,a_0+\frac 12)$, the trajectory 
$\{\varphi_t\}_{t\geq 0}$ intersects the semicircle $\{z\mid |z-a_0|=1\}$; let
$z'$ be the point of intersection and 
let $\xi\in (-1,1)$ be such that $a_0+\xi$ is real part of $z'$. 
If $w>a_0$, which is in particular necessarily 
the case if $a_0=2$, then we have $\xi \leq \frac 12 < \mu$ and hence 
$z'\in \{a_0+z\mid z\in C\}$. 
We may therefore suppose that $a_0\geq 3$ and $w<a_0$. We next observe
that as $u\leq \lambda$ and $w>2$, if $\xi_0$ is such that $a_0+\xi_0$ is the 
point of intersection of the geodesic joining $\lambda$ and $a_0-\frac 12$ 
then $\xi_0<\xi<0$. A direct computation shows
that $\xi_0=-\mu$ if $a_0=3$ and if $a_0>3$  in fact $\xi_0>-\mu$. 
This shows that $z'$ is contained in $\{a_0+z\mid z\in C\}$. Let $\zeta'$ 
be the unit tangent vector at $z'$, tangent to $\varphi$. 

Now let $z_1=-1/(z'-a_0)$, and $\zeta_1$ be the tangent vector at $z_1$,
which is the image of the  $\zeta'$ under the
map $z\mapsto  -1/(z-a_0)$. Then $z_1\in C$, and we claim that $(z_1,\zeta_1)
\in \Phi$, namely that $\tilde \varphi (z_1,\zeta_1)$ is $H$-reduced.  
Clearly the repelling and attracting endpoints of $\tilde \varphi 
(z_1,\zeta_1)$ are given by $u_1=-1/(u-a_0)$ and $w_1=-1/(w-a_0)$ respectively. 
Since $|w-a_0|<\frac 12$ we have $|w_1|=|1/(w-a_0)|>2$, as required. Since 
$a_0\geq 2$ and $u\leq \lambda <1$,  $u_1=-1/(u-a_0)$ is positive, and we only 
need to confirm that the appropriate upper bounds hold for $u_1$, depending 
on the sign of $w_1$. 
If $w>a_0$ then  $w_1<0$, and we have $u_1=\frac 1{a_0-u}
\leq \frac 1{2-\lambda}=1- \lambda $, and on the 
other hand if $w<a_0$ then  $w_1>0$ and $a_0\geq 3$, and in this 
case we have 
$u_1=\frac 1{a_0-u}\leq \frac 1{3-\lambda}=\lambda$. Thus 
$(z_1, \zeta_1)\in \Phi$, which proves the claim.  

Let $t_1>0$ be such that $\varphi_{t_1}=z'$, the latter being as above.
We note also that for $w_1$, which is the attracting endpoint of 
of $\tilde \varphi (z_1, \zeta_1)$, the negative Hurwitz continued fraction
expansion is given 
by $[a_1, a_2, \dots]$. An analogous arguments works by symmetry with 
when $a_0<0$ and we get $(z_1,\zeta_1)\in \Phi$, and $t_1>0$ such that 
$z_1=-1/\varphi_{t_1}$, $\zeta_1$ is the image of the unit tangent to 
$\varphi$ at $\varphi_{t_1}$ under the corresponding map, and the 
attracting endpoint of $\tilde \varphi (z_1, \zeta_1)$, is given 
by $[a_1, a_2, \dots]$.

Repeating the procedure we get a sequence 
$\{(z_j,\zeta_j)\}$ in $\Phi$ 
and a sequence $\{t_j\}$ of positive numbers such that
$z_j=-1/(\varphi_{t_j}(z_{j-1})-a_{j-1})$ 
and the tangent to $\varphi$ at 
$\varphi_{t_j}(z_{j-1})$ is mapped to $\zeta_j$ under the corresponding
tangent map (in particular the pair $\varphi_{t_j}(z_{j-1})$ together with 
the unit tangent to $\varphi $ at the point  
is equivalent to $(z_j, \zeta_j)$ under the $\Gamma$ action).  Also, 
the attracting endpoints of $\tilde \varphi (z_j,\zeta_j)$ have the 
continued fraction expansion $[a_j, a_{j+1}, \dots ]$, for all $j$. 

We call $\{t_j\}$ as above the sequence of {\it return times} corresponding
to the $(z_0,\zeta_0)$ or equivalently to the reduced geodesic $\varphi= 
\tilde \varphi (z_0,\zeta_0)$. 

Now let  $(z_0,\zeta_0)
\in \Phi$,  $\varphi =\tilde \varphi (z_0,\zeta_0)$, and $\psi$ the image
of $\varphi$ under the quotient map 
$\eta:\H^2\to M$. Let $\varphi 
=\{\varphi_t\}_{t\in \R}$; we equip $\psi$ with the parametrization given by
$\psi_t=\eta (\varphi_t)$ for all $t\in \R$. In the following we shall 
concern ourselves only with the (forward) trajectory $\{\psi_t\}_{t\geq 0}$. 
Using 
the return times of $\varphi$ we divide the trajectory into parts $\psi^{(j)}$,
$j\geq 0$; the indexing starting with $0$ is chosen for notational 
convenience with respect to the related indices, as will be seen below. 
For each $j\geq 0$ let  
$$\varphi^{(j)}=\{\tilde \varphi (z_{j}, \zeta_{j})(t)\mid 0
\leq t < t_{j+1}\} \ \  \mbox{ \rm and } \ \  \psi^{(j)}=\eta (\varphi^{(j)}). $$

In the following propositions we collect properties of the segments 
$\varphi^{(j)}$ to be used in our counting results in the next section. 
In the following let $[a_0, a_1, \dots]$ be the negative Hurwitz continued fraction 
expansion of the attracting endpoint of $\tilde \varphi (z_0,\zeta_0)$. 
 Also let $u_j$ 
and $w_j$ be the repelling and attracting endpoints of $\tilde \varphi 
(z_j,\zeta_j)$.

\begin{proposition}\label{intersect} 
Let $\delta \in (0,1)$ and $j\geq 0$ be given. If $\varphi^{(j)} 
\cap \H_\delta\neq \emptyset$ then $|a_j|>2\delta^{-2}
 +\lambda-\frac 32$, and if  $|a_j|>2\delta^{-2} +\lambda+\frac 12$ 
then $\varphi^{(j)} \cap \H_\delta\neq \emptyset$. When 
nonempty the intersection is an arc along the geodesic. 

\end{proposition}

\proof We recall that $\varphi^{(j)}$ is the segment of the (euclidean) 
semicircle joining 
$u_j$ and $w_j$ lying between $C$ and $a_j+C$ (both of which it intersects),
and its intersection with $\H_\delta$ is the same as that of the semicircle;
in particular, when nonempty it is an arc.   
Since $|w_j-a_j|<\frac 12$ and sgn$(w) u_j\in [\lambda -1, \lambda]$ it follows 
that the radius of the circle is between $\frac 12(|a_j|-\frac 12 
-\lambda)$ and $\frac 12(|a_j|+\frac 12 -(\lambda-1))$. The assertion in 
the proposition is immediate from these observations. \qed

\begin{remark}\label{rem2}
{\rm In particular Proposition~\ref{intersect} shows that $\{\psi\}_{t\geq 0}$ 
is bounded (has compact closure) in $M$ if and only if the sequence of 
partial quotients $\{a_j\}$ is bounded. Together with Remark~\ref{rem1} 
this shows that if $g\in G$
is such that $g(0)$ and $g(\infty)$ are irrational numbers whose minus
continued fraction expansion with respect to the Hurwitz algorithm 
have bounded  partial quotients  then for $Q=Q_g$ the inequality 
$|Q(p)|<\epsilon$ has no nonzero solution for sufficiently small 
$\epsilon >0$. We note also that this also shows, independently, that 
the partial quotients in the Hurwitz expansion of an irrational number are 
bounded if any only if the number is badly approximable in the usual sense. 
}
\end{remark}

\begin{proposition}\label{prop:mdelta}
Let $0<\delta <1$. Then for any $j\geq 0$ the arc $\varphi^{(j)}
\cap \H_\delta$ (when 
nonempty) is a connected component of $\eta^{-1} (\psi^{(j)}\cap M_\delta)$. 
Moreover,
if $\delta < (1-\mu^2)^{1/4}$ then it is the only connected component of 
$\eta^{-1} (\psi^{(j)}\cap M_\delta)$ contained in $\H_\delta$.
\end{proposition}

\proof For $\delta <1$, for any $\gamma\in \Gamma$ which does not
 leave $\H_\delta$ invariant
$\gamma \H_\delta$ is contained in $\{z=x+iy\in \H^2 \mid y<1\}$. Hence 
the intersection 
of $ \varphi^{(j)}$ with any such $\gamma \H_\delta$ is separated from 
$\varphi^{(j)}\cap \H_\delta$. This proves the first statement. If 
$\delta < (1-\mu^2)^{1/4}$
then for any $\gamma\in \Gamma$ which does not leave $\H_\delta$ invariant
$\gamma \H_\delta$ is contained in $\{z=x+iy\in \H^2 \mid y< \sqrt {1-\mu^2}\}$.
Since the endpoints of $\varphi^{(i)}$ have the $y$ coordinates greater than 
$\sqrt {1-\mu^2}$, $\varphi^{(i)}$ does not intersect $\gamma \H_\delta$ which
is different from $\H_\delta$. This completes the proof. \qed

\medskip
We shall denote by $d(z, z')$, where $z,z'\in \H^2$, the (hyperbolic) distance  
between $z$ and $z'$ in $\H^2$.  Let $\{\chi_j\}$ be the sequence of numbers
defined by $\chi_j=\frac 12 \log 3\sqrt 5$ if $|a_j|\geq 3$ and
$\chi_j =\log 2-\frac 12  d\left(\mu+i\sqrt {1-\mu^2},\frac 12(3+i\sqrt 3)\right)$
if $|a_j|=2$. 

\begin{proposition}\label{length}
$ - 2\chi_j\leq t_j-2\log |a_j| \leq  \log 3\sqrt 5+\log 
\left (\frac 34 + \sqrt {\frac 12}\right )$, for all $j\geq 0$. 
\end{proposition}

\proof Recall that each $\varphi^{(j)}$ is a geodesic segment joining a point of 
$C$ to another on $a_j+C$, and $t_j$ is its length. The distances of 
the initial point and the end point of the segment from $i$ and $a_j+i$ 
respectively are bounded by $d(\mu+\sqrt {1-\mu^2},i)$, 
where, as before, $\mu =\frac{23-3\sqrt 5}{22}$ is the $x$ co-ordinate 
of the right endpoint of $C$. The latter distance 
is $\frac 12\log \frac{1+\mu}{1-\mu}$, and a numerical calculation shows 
that it equals $\frac 12 \log 3\sqrt 5$.  Now let $d_j=d(i,a_j+i)$.  
Then from the above observations 
we get that
\begin{align*}
 d_j-\log 3\sqrt 5 \leq  t_j\leq d_j+\log 3\sqrt 5 ,
\end{align*}
 
for all $j$. By a 
standard formula for distances in $\H^2$ (see \cite{K}) we have
\begin{align*}
d_j&= \log (\frac 12 |a_j^2|+\frac 12|a_j|{\sqrt {|a_j|^2+4}+1})\\
~&= 2\log |a_j| +\log (\frac 12+\frac 12\sqrt{1+4|a_j|^{-2}}+|a_j|^{-2}) 
\end{align*}
\noindent and since $|a_j|\geq 2$ for all $j$ the second term is 
(positive and) at 
most $\log \left(\frac 34 + \sqrt {\frac 12}\right)$. 
Combining, we get that
\begin{align*}
 -\log 3\sqrt 5\leq t_j- 2\log |a_j|\leq   
\log 3\sqrt 5+\log 
\left(\frac 34 + \sqrt {\frac 12}\right),
\end{align*}
for all $j\geq 0$. From the definition 
of $\chi_j$'s we see that this proves the 
proposition in the case $|a_j|\geq 3$, and also the second inequality 
when $|a_j|=2$. 

It remains to prove the first inequality in the case when $|a_j|=2$. 
For this we note that $\varphi^{(j)}$ is a segment joining a point of $C$
to a point of $\pm 2+C'$, where $ C'=\{z=x+iy\in C\mid |x|\leq \frac 12\}$.
It can be verified directly that the geodesic joining the endpoints
$\mu +i\sqrt {1-\mu^2}$ and  $\frac 12(3+i\sqrt 3)$ of $C$ and $2+C'$ is, 
along with its mirror image, is the shortest of the segments as above. 
Thus $ t_j\geq d\left(\mu+i\sqrt {1-\mu^2},\frac 12(3+i\sqrt 3)\right).$ 
Hence 
\begin{equation*}
 t_j-2\log |a_j|\geq d\left(\mu+i\sqrt {1-\mu^2},\frac 12(3+i\sqrt 3)
\right)-2\log 2 =-2\chi_j,
\end{equation*}
by the definition of $\chi_j$ when $|a_j|=2$. 
\qed     

\begin{remark}\label{chi2}
{\rm Let $\chi =\log 2 -\frac 12 d\left(\mu+i\sqrt {1-\mu^2},
\frac 12(3+i\sqrt 3)\right)$, namely the value of $\chi_j$ when $|a_j|=2$. 
We note that $
 \frac 12 \log \frac {16}{11}\leq \chi \leq \frac 12 \log 
\frac 32$. The semicircular 
geodesic segment joining $\frac 12(3+i\sqrt 3)$ to $\frac 34+i\frac 
{\sqrt 7}4$ may be seen to 
be orthogonal to the unit circle at the latter point. Hence we have
\begin{equation*}
  d\left(\mu+i\sqrt {1-\mu^2},\frac 12(3+i\sqrt 3)\right)\geq d\left(\frac 34+i\frac 
{\sqrt 7}4,\frac 12(3+i\sqrt 3)\right).
\end{equation*}
By a direct computation we see that $
 d\left(\frac 34+i\frac {\sqrt 7}4,
\frac 12(3+i\sqrt 3)\right)=\log \frac 
{2+\sqrt 7}{\sqrt 3}\geq \log \frac 83$.
Hence 
$ \chi \leq \log 2 -\frac 12
\log \frac 83 = \frac 12 \log \frac 32$.
Also, 
\begin{align*}
d\left(\frac 34+i\frac {\sqrt 7}4, \mu+i\sqrt {1-\mu^2}\right) &=
d\left(\frac 34+i\frac {\sqrt 7}4,i\right)-d\left(\mu+i\sqrt {1-\mu^2}, i\right)\\\\
&=\log \sqrt 7 -\frac 12\log 3\sqrt 5 ,
\end{align*}
from which we get that
 
 \begin{align*}
  d\left(\mu+i\sqrt {1-\mu^2},\frac 12(3+i\sqrt 3)\right)
  &\leq \log  \frac {2+\sqrt 7}{\sqrt 3} +\log \sqrt 7-
\log\  ^4\!\!\!\sqrt 45 \leq \log \frac {11}4,
 \end{align*}
as may be directly verified. Thus we have 
$ \chi \geq \log 2-\frac 12 \log  \frac {11}4 =\frac 12 \log \frac {16}{11}$.

}
\end{remark}

\section{Solutions of $H$-reduced quadratic inequalities} 

We say that an element $g\in G$ is $H$-reduced if the geodesic $\{ga_tK\}$
(under the identification as in \S~3) is $H$-reduced. Now let 
$g\in G$ be such that $g(\infty)$ is an irrational number and let 
$Q=Q_g$; such a quadratic form may be called an {\it $H$-reduced quadratic 
form}. Let $[a_0, a_1, \dots]$ be the minus continued fraction expansion 
of $g(\infty)$ with respect to the Hurwitz algorithm. For all $n\in \N$ let 
$$\alpha_n=\sum_{j=0}^{n-1}\log |a_j|.$$ 
Also for $\delta >0$ and $n\in \N$ let 
\begin{align*}
~~~~e(\delta,n)= \# \{0\leq j\leq n-1\mid |a_j|>2\delta^{-2}+\frac 12 
+\lambda\} \mbox{ \rm and } \\
f(\delta,n)= 
\#\{0\leq j\leq n-1\mid |a_j|>2\delta^{-2}-\frac 32 +
\lambda\},~~~~~
\end{align*}
  where, as in \S~3, $\lambda =\frac 12(3-\sqrt 5)$. Also let $c_0=
\frac 12\log 3\sqrt 5+\frac 12\log \left (\frac 34 + \sqrt {\frac 12}
\right )$.

\begin{theorem}\label{thm:reduced}
Let $0< \delta <\sqrt {\frac 2\pi}$ and $\kappa >0$ be given, and for 
$\rho>0$ let  $$G(\rho)= \{p\in \mathfrak p\mid 0<|Q(p)|
<\frac 12 \delta^2, L_g^-(p) >\kappa 
\mbox { \rm and }  \norm p \norm \leq \rho\}.$$
There exist constants $\theta\geq 0$ and $\nu\in \N$ such that  
the following holds:
 
i) if  $\log \rho \geq \alpha_n+c_0n+\theta $ 
then $\# G(\rho)\geq e(\delta, n)-\nu$, and 

ii) if $\log \rho \leq \alpha_n-\sum_{j=0}^{n-1}\chi_j
-\theta$ then $\# G(\rho)\leq f(\delta, n)+\nu$.
\end{theorem}

\proof Let $\delta >0$ be given. By Corollary~\ref{cor3}, for any 
$\rho >0$ the cardinality of $G( \rho)$ 
differs from 
$c_g^+(2\log \rho)$ by a bounded amount, where $c_g^+(2\log \rho)$ stands for 
the number of connected components of the set 
$\{t\in [0,2\log \rho]\mid ga_tB_\delta\cap \Z^2\neq \emptyset\}$. 
 Hence it suffices
to show that the estimates as in assertions~(i) and~(ii) hold for
$c^+_g(2\log \rho)$  (for the given $\delta$, not included in the notation) 
in place of $\# G(\rho)$, with a $\nu'$ in place of $\nu$.

We note that  $t\geq 0$, $ga_tB_\delta\cap \Z^2\neq \emptyset$ 
if and only if $a_t^{-1}g^{-1}\gamma e_1 \in B_\delta$ which, as seen in \S~3, is 
equivalent to $\eta 
(ga_t(i))\in M_\delta$. 

Now let $(z,\zeta) \in \Phi$ be such that $\tilde \varphi (z,\zeta)$ 
is equivalent to the geodesic joining $g(0)$ and $g(\infty)$. Let $\varphi=
\{\varphi_t\}=\tilde \varphi(z,\zeta)$. Let $\{t_j\}_{j=0}^\infty$ denote the 
corresponding
sequence of return times (see \S~3). Also let $t'\in \R$ be such that
 $g(i)=\varphi_{t'}$; such a $t'$ exists since $g(i)$ is contained in the 
geodesic joining $g(0)$ and $g(\infty)$. We choose $\theta =\frac{\disp{|t'|}}{\disp 2}$ and 
$\nu'$ to be the infimum of $j$ such that $t_j\geq t'$. We show that 
assertions~(i) and~(ii) hold for these choices. 

Now let $\rho$ be as in~(i), with $n\in \N$.  Suppose
first that $t'<0$. Then  $$2\log \rho \geq  2\alpha_n +2c_0n-t'\geq 
\sum_{j=0}^{n-1} (2\log |a_j|+2c_0)-t'> \sum_{j=0}^{n-1} t_j-t'$$ and hence 
$c_g^+(2\log \rho)$ contains $\{t\in [t_j, t_{j+1})\mid \varphi_t\in \H_\delta\}$, 
for all $j=0, \dots , n-1$. By 
Proposition~\ref{prop:mdelta}, this implies  $c_g^+(2\log \rho)$ is at least as 
much as the number of $j\leq n-1$ for which
 $\varphi^{(j)}$ intersects $\H_\delta$. By Proposition~\ref{intersect} this 
number is at least $e(\delta, n)$, which proves the claim
in the case at hand. Now suppose that $t'\geq 0$. In this 
case $c_g^+(2\log \rho)$ is seen to be 
at least the number
of those $j$'s for which $t_j>t'$ and $\varphi^{(j)}$ intersects $\H_\delta$.
 By Proposition~\ref{intersect} the 
number of $j\leq n-1$ for which $\varphi^{(j)}$ intersects $\H_\delta$ is 
at least $e(\delta, n)$. This shows that $c_g^+(2\log \rho)$ is  
at least  $e(\delta, n)-\nu$ which, as  noted above, 
proves (i). 

Next let $\rho$ be as in~(ii) with $n \in \N$. Then $$2\log \rho \leq  
2\alpha_n -2\sum_{j=0}^{n-1}\chi_j-t'\leq \sum_{j=0}^{n-1} (2\log |a_j|-2\chi_j)-t'< 
\sum_{j=0}^{n-1} t_j-t'.$$
Since $\sqrt{\frac 2\pi}< (1-\mu^2)^{1/4}$, by Proposition~\ref{prop:mdelta}
this implies that  $c_g^+(2\log \rho)$ is at most the number of
 $j\leq n-1$ for which
 $\varphi^{(j)}$ intersects $\H_\delta$, and  Proposition~\ref{intersect} it 
 is at most $f(\delta, n)$. This proves~(ii). \qed

\section{Solutions of quadratic inequalities - the general case}

Consider a binary quadratic form $Q(x,y)$, which is nondegenerate and not
a scalar multiple of a form with rational coefficients. Then upto a 
scalar multiple $Q$ is given by $Q(x,y)= (ay+bx)(cy+dx)$, 
for all $x,y\in \R$, where $ad-bc=1$, $b\neq 0$  and $\frac ab $ is irrational. 
We shall therefore consider only forms $Q$ satisfying these conditions 
on the coefficients. Thus 
 $Q(x,y)= Q_g(xe_1+ye_2)$, where $g  = \left(\begin{matrix}
~-c & a \\
d & ~-b \\
\end{matrix}\right)$
is an element in $G$.  We note that $g(\infty)=-\frac ab$, which by hypothesis
is an irrational number. The element $g$ may not be $H$-reduced, namely 
$\{ga_tK\}$ may not be a
$H$-reduced geodesic.  However  it is equivalent under the action of $\Gamma$ 
to an $H$-reduced geodesic. Thus there exists a $\gamma \in \Gamma$
such that $\{\gamma g a_tK\}$ is $H$-reduced. Let $g'=\gamma g$ and $Q'=
Q_{g'}$; then $g'$ is $H$-reduced, and we shall call $Q'$ a {\it reduced version 
of $Q$} and $\gamma$ an {\it $H$-reducing element} for $Q$.  

We note that the factors $L_g^+$ and $L_g^-$ of $Q_g$ as introduced earlier 
are now given by $L_g^+(xe_1+ye_2)=ay+bx$ and  $L_g^-(xe_1+ye_2)
=cy+dx$, for all $x,y\in \R$. Let 
$L_{g'}^+$ and $L_{g'}^-$ be the corresponding linear forms for $g'$. 
 
\medskip
Let  $[a_0,a_1, \dots, ]$ and  $[a_0',a_1', \dots, ]$
be the negative Hurwitz continued fraction expansions of $\frac ab$ (which is 
irrational) and $g'(\infty)$ respectively. 
Since  $g(\infty)
=-\frac ab $, it follows that $[-a_0',-a_1', \dots, ]$ is the 
continued fraction expansion of $g'(\infty)$. Since $g'(\infty)=\gamma g
(\infty)$ we get also that there exists $m\in \Z$ such that 
$a'_j=-a_{j+m}$ for all large $j$. Now let 
$$\alpha^+=\limsup \frac 1n \sum_{j=0}^{n-1} \log |a_j|=\limsup \frac 1n 
\sum_{j=0}^{n-1} \log |a'_j| $$
and 
$$\alpha^-=\liminf \frac 1n \sum_{j=0}^{n-1} \log |a_j|=\liminf \frac 1n 
\sum_{j=0}^{n-1} \log |a'_j|. $$
Let $\{\chi_j\}$ be the sequence, as before (defined by $\chi_j=
\frac 12 \log 3\sqrt 5$ if $|a_j|> 2$ and $\chi_j =\log 2-\frac 12  
d(\mu+i\sqrt {1-\mu^2},\frac 12(3+\sqrt 3))$ if $|a_j|=2$) and let   
$$\omega_n = \sum_{j=0}^{n-1} (\log |a_j|-\chi_j) \mbox{ \rm for all }n
\in \N, \  \mbox{ \rm  and } \omega =\liminf \frac 1n\omega_n . $$
It may be borne in mind that each of $\alpha^+, \alpha^-$ and $\omega$ 
can be infinite. On the other hand $\alpha^+\geq \alpha^-\geq \log 2$. 
Also, as 
 $$\log 2 -\chi \geq \log 2 -\frac 12 \log \frac 32 > \log 3 - \frac 12 
\log 3 \sqrt 5
=\frac 14 \log \frac 95,$$  we have $\log |a_j| -\chi_j\geq \frac 14 
\log \frac 95$ for all $j$, and hence $\omega \geq \frac 14 
\log \frac 95$. We note also that similarly $\omega \geq \eta \alpha^-$, where 
$\eta =(\log \frac 95)/(4 \log 3)>\frac 18$. 

Also, for any $\delta>0$
let $$e^-(\delta)={\displaystyle \liminf_{n\to \infty}}\, \frac 1n e(\delta, n),  
 \quad e^+(\delta)={\displaystyle \limsup_{n\to \infty}}\, \frac 1n e(\delta, n)$$
 \  $$f^-
(\delta)= {\displaystyle \liminf_{n\to \infty}}\,
\frac 1n f(\delta, n) \mbox{ \rm and }   
 f^+(\delta)={\displaystyle \limsup_{n\to \infty}}\, \frac 1n f(\delta, n).$$
We note that $\alpha^+$, 
$\alpha^-$, $\omega$,  $e^+(\delta)$ and $e^-(\delta)$ are determined by 
$\{a_j\}$, and hence by $\frac ab$. 
 
\begin{corollary}\label{asymptotic}
Let $Q$ be a binary quadratic form, as above, given by $Q(x,y)= 
(ay+bx)(cy+dx)$, for all $x,y\in \R$, where $ad-bc=1$,  $b\neq 0$ 
and $\frac ab $ is irrational. Let the notation
$\alpha^+$, $\omega$, $e^+(\delta)$ and $e^-(\delta)$ be 
as above, and let $c_0=\frac 12\left(\log 3\sqrt 5+\log (\frac 34 + 
\sqrt {\frac 12})\right)$, as before. 
Let $0<\delta <\sqrt{\frac 2\pi}$ and $\kappa>0$ be given, and for 
any $\rho$ let 
$$G(\rho)= \{p=xe_1+ye_2\in \mathfrak p\mid 0<|Q(p)|<\frac 12 \delta^2,\  cy+dx  
>\kappa \mbox { \rm and }  \norm p \norm \leq \rho\}. $$ Let 
$\epsilon >0$ be arbitrary. Then we have the following: 

i) if $\alpha^+<\infty$ then there exists $\rho_0$ such that 
for all $\rho \geq \rho_0$ we have  
$$\# G(\rho) \geq (1-\epsilon) \frac {e^-(\delta)}{(\alpha^++c_0)}\log \rho
;$$  

ii) given $0 <M\leq \omega $ there exists $\rho_0$ such that 
for all $\rho \geq \rho_0$ we have  
$$\# G(\rho) \leq (1+\epsilon) \frac {f^+(\delta)+\epsilon }
{M}\log \rho.$$
In particular this holds for all $0<M\leq \eta \alpha^-$. 
\end{corollary}

\proof 
Let $\gamma\in \Gamma$ be an $H$-reducing element for $Q$, $g'=\gamma g$ 
and $Q'=Q_{g'}$. 
For any $\rho >0$ let $G'(\rho)= \{p\in \mathfrak p\mid 0<Q'(p)<\delta^2, 
L_{g'}^-(p) >\kappa \mbox { \rm and }  \norm p \norm \leq \rho\}$. Then 
we have $Q(p)=Q'(\gamma p)$ and $L_g^-(p)=L_{g'}^-(\gamma p)$ for all 
$p\in \Z^2$,  and hence  
$ \# G'(\norm \gamma \norm^{-1} \rho)\leq \# G(\rho)\leq \# G'(\norm 
\gamma \norm \rho)$ for all $\rho >0$, where 
$\norm \gamma \norm $ denotes the operator norm of $\gamma$ (as an 
element of $G$).  This shows that it suffices 
to prove the assertions in the corollary with $G'(\rho)$ in place of 
$G(\rho)$, namely in the case when $\gamma$ 
is the identity element. In other words we may assume, as we shall, 
that $g$ is $H$-reduced. 

Let $\theta $ and $\nu$ be as in Theorem~\ref{thm:reduced}. Let $\upsilon
=\sqrt {1+\epsilon}$.  There exists
$n_0$ such that for all $n\geq n_0$ we have $e(\delta,n)-\nu\geq 
{\upsilon^{-1}}e^-(\delta)n$ and $\alpha_{n+1}+c_0(n+1)+\theta \leq 
 \upsilon (\alpha^++c_0)n$. Let $\rho_0>0$ be such that $\log \rho_0
=\alpha_{n_0}+c_0n_0+\theta$. Consider $\rho \geq \rho_0$. Then there exists 
$n\geq n_0$ such that $\alpha_n+c_0n+\theta \leq \log \rho \leq \alpha_{n+1}
+c_0(n+1)+\theta$. Then by Theorem~\ref{thm:reduced} we have $\# G(\rho)
\geq e(\delta , n)-\nu \geq  {\upsilon^{-1}}e^-(\delta)n $. 
Also by choice $\log \rho \leq \alpha_{n+1} +c_0(n+1)+\theta\leq 
 \upsilon (\alpha^++c_0)n$, and hence $n\geq \log \rho/ 
\upsilon (\alpha^++c_0)$. 
Thus we get that $$\# G(\rho)\geq {\displaystyle {\upsilon^{-2}}
\frac {e^-(\delta)}
{(\alpha^++c_0)}\log \rho \geq (1-\epsilon)\frac {e^-(\delta)}
{(\alpha^++c_0)}\log \rho},$$ which proves~(i). 

Next we choose $n_0$ such that for all $n\geq n_0$ we have 
$f(\delta,n)+\nu\leq 
 {\upsilon}(f^+(\delta)+\epsilon)n$ and $\omega_{n-1}
-\theta \geq  \upsilon^{-1} Mn$.  Let $\rho_0>0$ be such that $\log \rho_0
=\omega_{n_0}-\theta$, and consider $\rho\geq \rho_0$. Since $\omega>0$, 
$\omega_n\to \infty$ and hence there exists $n\geq n_0+1$ such that 
$\log \rho \leq \omega_n-\theta$; we pick the least integer $n\geq n_0+1$  
with this 
property, so $\log \rho \geq \omega_{n-1}-\theta$.  
By Theorem~\ref{thm:reduced} we have $\# G(\rho)
\leq f(\delta , n)+\nu \leq  {\upsilon}(f^+(\delta)+\epsilon)n $. 
Also since $\log \rho \geq \omega_{n-1}-\theta\geq 
{\upsilon^{-1}} Mn$, we get that $n\leq \frac  
\upsilon M\log \rho$. 
This yields $$\# G(\rho)\leq {\displaystyle {\upsilon^2}\frac 
{f^+(\delta)+\epsilon}
{M}\log \rho = (1+\epsilon)\frac {f^+(\delta)+\epsilon}
{M}\log
 \rho},$$ which proves~(ii). \qed

\bigskip
\noindent{\it Proof of Theorem~\ref{thm:main}}: 
Note that the values of $\alpha^+$, $\alpha^-$,
$e(\delta)$ and $f(\delta)$ involved in the statement of the theorem,
coincide for the Hurwitz and negative Hurwitz expansion.
%
The theorem now follows from
Corollary~\ref{asymptotic} when we interchange the role of $x$ and $y$ and 
replace $\delta $ by $\sqrt {2\delta}$; the constants $e(\delta)$ and 
$f(\delta)$ as in Theorem~\ref{thm:main} are respectively at most and 
at least as much as the corresponding constants in  
Corollary~\ref{asymptotic}, since $\lambda <\frac 12$.

\qed

\medskip

The following special case may be worth emphasizing, on account of its 
comparability with the fact that if $[a_0,a_1, \dots, ]$ is bounded then 
for sufficiently small $\delta >0$ the set of solutions $G(\rho)$ as above 
is   empty. 

\begin{corollary}\label{zerodensity}
Let the notation be as above. If  $f^+(\delta)=0$ then 
$$\lim_{\rho \to \infty}\frac{\#G(\rho)}{\log \rho} =0;$$
in particular,  if $S$ is
a subset of $\N$ with zero 
upper density and $\{a_j\}$ is bounded on  the complement of $S$,  then  
this conclusion holds, for 
all sufficiently small $\delta>0$. 
\end{corollary}

\proof The proof is immediate from Corollary~\ref{asymptotic}. \qed

%


\smallskip
\bibliography{ref}
\bibliographystyle{plain}







\bigskip
\begin{flushleft}
Manoj Choudhuri  \hfill S.G. Dani

Centre for Applicable Mathematics \hfill Department of Mathematics

Tata Institute of Fundamental Research \hfill Indian Institute of Technology Bombay

Yelahanka, Bangalore 560065, India. \hfill Powai, Mumbai 400076, India.

manoj@math.tifrbng.res.in \hfill sdani@math.iitb.ac.in
\end{flushleft}
\end{document}